\documentclass[submission,copyright,creativecommons]{eptcs}
\providecommand{\event}{GASCom 2024} 

\usepackage{iftex}

\ifpdf
  \usepackage{underscore}         
  \usepackage[T1]{fontenc}        
\else
  \usepackage{breakurl}           
\fi

\usepackage{graphicx}
\usepackage{amsmath}
\usepackage{amssymb}
\usepackage{tikz-cd}

\newcommand{\PLR}{\mathcal{P}^{LR}_k(\Omega)}
\newcommand{\Zd}{\mathbb{Z}[\delta]}
\newcommand{\PLRn}{\widehat{\mathcal{P}}^{LR}_{n+2}(\Omega)}
\newcommand{\PLRk}{\widehat{\mathcal{P}}^{LR}_{k}(\Omega)}
\newcommand{\TLR}{T^{LR}_k(\Omega)}

\newcommand{\locb}{\mathcal{L}_\bullet^{\circ}}

\newcommand{\TLD}{\mathrm{TL}(\widetilde{D}_{n+2})}

\newcommand{\Dn}{\widetilde{D}_{n+2}}
\newcommand{\DD}{\mathbb{D}(\widetilde{D}_{n+2})}
\newcommand{\dbD}{\mathtt{ad}(\widetilde{D}_{n+2})}

\newcommand{\FC}{\operatorname{FC}}

\newtheorem{Theorem}{Theorem}[section]

\newtheorem{Definition}[Theorem]{Definition}

\title{Diagram Calculus for the Affine Temperley--Lieb \\Algebra of Type $D$}
\author{Riccardo Biagioli
\institute{Università di Bologna\\ Bologna, Italy}
\email{riccardo.biagioli2@unibo.it}
\and
 Giuliana Fatabbi
\institute{Università degli Studi di Perugia\\
Perugia, Italy}
\email{giuliana.fatabbi@unipg.it}
\and
 Elisa Sasso
\institute{Università di Bologna\\
Bologna, Italy}
\email{elisa.sasso2@unibo.it}
}

\def\titlerunning{Diagram Calculus for the Affine Temperley--Lieb Algebra of Type $D$}
\def\authorrunning{R. Biagioli, G. Fatabbi \& E. Sasso}
\begin{document}
\maketitle

\begin{abstract}
Let $(W,S)$ be a Coxeter system of affine type $\widetilde{D}$, and let $\mathrm{TL}(W)$ the corresponding generalized Temperley-Lieb algebra. In this extended abstract we define an infinite dimensional associative algebra made of decorated diagrams which is isomorphic to $\mathrm{TL}(W)$. Moreover, we describe an explicit basis for such an algebra of diagrams which is in bijective correspondence with the classical monomial basis of $\mathrm{TL}(W)$, indexed by the fully commutative elements of $W$.
\end{abstract}

\section{Introduction}


The Temperley-Lieb algebra is a very classical mathematical object studied in algebra, combinatorics, statistical mechanics and mathematical physics, introduced by Temperley and Lieb in 1971 \cite{TemperleyLieb}. Thanks to Kauffman \cite{Kauffman} and Penrose \cite{Penrose}, it was showed that the Temperley-Lieb algebra can be realized as a diagram algebra, that is an associative algebra with a basis given by certain diagrams on the plane. On the other hand,  Jones presented the Temperley-Lieb algebra in terms of abstract generators and relations. In \cite{Joneshecke}, he also showed that this algebra occurs naturally as a quotient of the Hecke algebra of type $A$. The realization of the Temperley-Lieb algebra as a Hecke algebra quotient was generalized by Graham in \cite{Graham}. He defined the so-called generalized Temperley Lieb algebra TL($\Gamma$) for any Coxeter system of type $\Gamma$ and showed that TL($\Gamma$) has a monomial basis indexed by the fully commutative elements of the underlying Coxeter group.
\begin{figure}[h]
    \centering
    \includegraphics[scale=0.7]{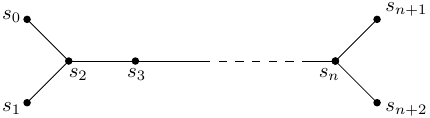}
    \caption{Coxeter graph of type $\Dn$.}
    \label{fig:graphD}
\end{figure}

During the years, diagrammatic representations for $\mathrm{TL}(\Gamma)$ have been found for each Coxeter system of finite type but only for two systems of affine type. More precisely, in \cite{Greengeneral}, \cite{GreenTLE} and \cite{Greencellular} Green defined a diagram calculus in finite Coxeter types $B$, $D$, $E$ and $H$. For affine types, in \cite{FanGreenAffine} Fan and Green provided a realization of TL($\widetilde{A}$) as a diagram algebra on a cylinder and, more recently, in \cite{ErnstDiagramI, ErnstDiagramII} Ernst represented TL($\widetilde{C}$) as an algebra of decorated diagrams. In this extended abstract, we present a new diagrammatic representation for $\mathrm{TL}(\Dn)$. Our method can be extended also to the affine case $\widetilde{B}$. Here we recall the presentation of TL($\Dn$) given by Green in \cite{GreenStar} that we consider as definition: TL($\Dn$) is the $\mathbb{Z}[\delta]$-algebra generated by $\{b_0, b_1, \ldots, b_{n+2}\}$ with defining relations:
	\begin{enumerate}
		\item[(d1)] $b_i^2=\delta b_i$ for all $i\in \{0,\ldots, n+2\}$,
        \item[(d2)] $b_i b_j=b_jb_i$ if $s_i$ and $s_j$ are not adjacent nodes in the Coxeter graph of type $\Dn$;
        \item[(d3)] $b_i b_j b_i=b_i$ if $s_i$ and $s_j$ are adjacent nodes in the Coxeter graph of type $\Dn$.
	\end{enumerate}

Similarly to Ernst in \cite{ErnstDiagramI, ErnstDiagramII}, we define the diagrams for our representation starting from the classical ones of type $A$ and adding decorations on the edges.


Then, as usual, we define a product of decorated diagrams by concatenation. This operation turns this set into an infinite dimensional $\mathbb{Z}[\delta]$-algebra, of which we consider a quotient modulo some new relations. Within this quotient, we consider a specific $\mathbb{Z}[\delta]$-subalgebra called \textit{algebra of admissible diagrams} and denoted by $\DD$. The main result of this extended abstract states that $\DD$ and $\TLD$ are isomorphic $\mathbb{Z}[\delta]$-algebras.

\section{Fully commutative elements of Coxeter groups}\label{fc}

Let $M$ be a square symmetric matrix indexed by a finite set $S$, satisfying $m_{ss}=1$ and, for $s\neq t$, $m_{st}=m_{ts}\in\{2,3,\ldots\}\cup\{\infty\}$. The \textit{Coxeter group} $W$ associated with the  \textit{Coxeter matrix} $M$ is defined by generators $S$ and relations $(st)^{m_{st}}=1$ if $m_{st}<\infty$. These relations can be rewritten more explicitly as $s^2=1$ for all $s$, and \[\underbrace{sts\cdots}_{m_{st}}  = \underbrace{tst\cdots}_{m_{st}},\]  where $m_{st}<\infty$, the latter being called \textit{braid relations}. When $m_{st}=2$, they are simply \textit{commutation relations} $st=ts$. 
For $w\in W$, the \textit{length} of $w$, denoted by $\ell(w)$, is the minimum length $l$ of an expression $s_1\cdots s_l$ of $w$ with $s_i\in S$. The expressions of length $\ell(w)$ are called \textit{reduced}.

\begin{Definition}
	\label{defi:FC}
	An element $w\in W$ is \textit{fully commutative} (FC) if any reduced expression of $w$ can be obtained from any other reduced expression of $w$ using only commutation relations.
\end{Definition}

\begin{figure}
    \centering
    \includegraphics[scale=0.6]{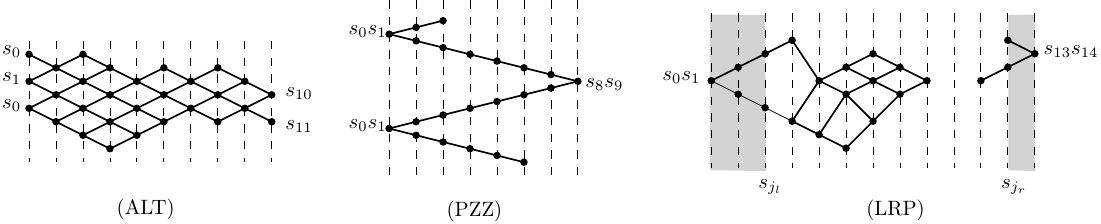}
    \caption{Some fully commutative heaps.}
    \label{fig:heaps}
\end{figure}

The concept of heap helps in studying problems related to full commutativity, for more details see for instance~\cite{BJNFC}. Briefly, given a reduced expression of $w=s_{i_1}\cdots s_{i_k}\in W$, its \textit{heap} is a poset on the index set $\{1,\ldots,k\}$ together with a labeling map. Heaps can be represented via Hasse diagrams; moreover, if $w\in \FC$, its heap does not depend on its reduced expression.
Fully commutative elements heaps of type $\Dn$ have been classified in  \cite[\S 3.2]{BJNFC}: they can be split in five disjoint families, depending on the shapes of their associated heaps, whose elements are respectively called Alternating elements (ALT), Left-Peaks (LP), Right-Peaks (RP), Left-Right-Peaks (LRP) and Pseudo-Zigzags (PZZ), see Figure \ref{fig:heaps}.

\section{Decorated diagrams}

A \textit{concrete pseudo k-diagram} consists of a finite number of disjoint plane curves, called \textit{edges}, embedded in a  box having $k$ nodes on the top (\textit{north}) face and $k$ nodes  on the bottom (\textit{south}) face. The nodes are  endpoints of  edges and  all other embedded edges must be closed (isotopic to circles) and disjoint from the box. We refer to a closed edge as a \textit{loop}. It follows that there cannot exist isolated nodes and  from each node  a single edge starts. By $\{a,b\}$ we mean an edge that joins the node $a$ to the node $b$.

\begin{figure}[hbtp]
	\centering
	\includegraphics[scale=0.45]{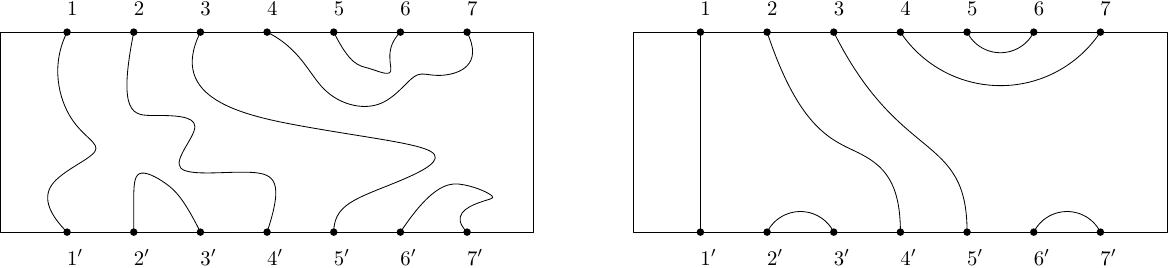}
	\caption{Two equivalent concrete pseudo 7-diagrams.}
	\label{kbox}
\end{figure}

We say that two concrete pseudo $k$-diagrams are \textit{(isotopically) equivalent} if one can be obtained from the other by isotopically deforming the edges such that any intermediate diagram is also a concrete pseudo $k$-diagram (Figure \ref{kbox}). We define a \textit{pseudo k-diagram} as an equivalence class of concrete pseudo $k$-diagrams with respect to isotopically equivalence. Given two of these diagrams $D,D'$, we define the product $D'D$ as the pseudo $k$-diagram obtained by placing $D'$ on top of $D$ so that node $i$ of $D$ coincides with node $i'$ of $D'$ and then rescaling.

Now let $D$ be a concrete pseudo $k$-diagram. Consider the set $\Omega=\{\bullet, \circ\}$ and the monoid $\Omega^*$. Our goal is to adorn the edges of $D$ with elements of $\Omega$ which we call \textit{decorations}. In particular, $\bullet$ is called a \textit{L-decoration} and $\circ$ is called \textit{R-decoration}. We call a $LR$-\textit{decorated pseudo k-diagram} a pseudo $k$-diagram decorated with these decorations up to certain rules that we do not list here. We denote the set of LR-decorated pseudo $k$-diagrams by $\TLR$ and define $\PLR$ to be the $\Zd$-module having the elements of $\TLR$ as a basis. 

\begin{figure}[hbtp]
	\centering
	\includegraphics[scale=0.7]{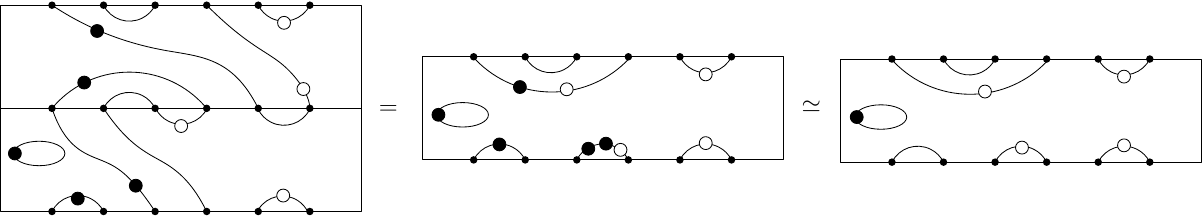}
	\caption{Product of two concrete decorated pseudo diagrams and its reduction.}
	\label{fig:product}
\end{figure}

As before, we define multiplication in $\PLR$ by concatenating two basis elements and then extend it bilinearly, see the first equality in Figure~\ref{fig:product}.
We can show that this product gives a structure of $\Zd$-algebra to $\PLR$, which is an infinite dimensional algebra.

\begin{figure}[hbtp]
	\centering
	\includegraphics[scale=0.6]{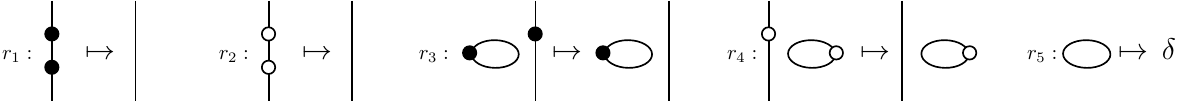}
	\caption{The defining relations of $\PLRk$.}
	\label{rel}
\end{figure}

Let $\PLRk$ be the $\Zd$-quotient algebra of $\PLR$ by the relations in Figure \ref{rel}.  We say that a LR-decorated diagram is \textit{irreducible} if there are no relations to apply. Similarly to \cite[Proposition 3.4.1]{ErnstDiagramI}, one can prove that the set of LR-decorated irreducible diagrams forms a basis for $\PLRk$. An example of irreducible diagram is in Figure~\ref{fig:product}, right.

\begin{figure}[h]
    \centering
    \includegraphics[scale=0.7]{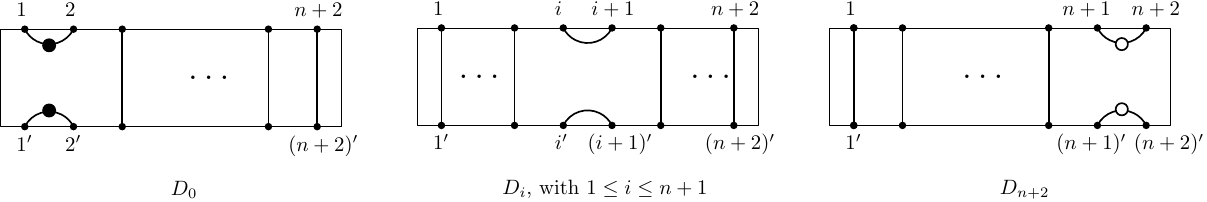}
    \caption{The simple diagrams.}
    \label{fig:simple}
\end{figure}

We are particularly interested in a special subset of irreducible diagrams, called the \textit{simple diagrams} $D_0, \ldots, D_{n+2}$, defined as in Figure \ref{fig:simple}. 
It is easy to prove that the simple diagrams satisfy the relations (d1)-(d3) of the $\TLD$ generators, defined in Section 1, simply replacing $b_i$ by $D_i$. Denote by $\DD$ the $\Zd$-subalgebra of $\PLRn$ generated as a unital algebra by the simple diagrams with multiplication inherited by $\PLRn$.\\

\subsection{Admissible diagrams}
We consider a subset of irreducible diagrams called \textit{admissible diagrams} and denoted by $\dbD$; here we state some fundamental properties the admissible diagrams must satisfy:
\begin{enumerate}
    \item the only loop edges that can occur are depicted in Figure \ref{loop-B};
    	\begin{figure}[h]
	\centering	\includegraphics[width=0.6\linewidth]{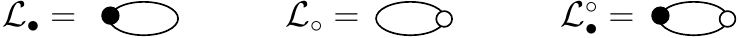}
	\caption{Allowable loops in $\tilde{D}$-admissible diagrams.}
	\label{loop-B}
\end{figure}

    \item the total number of $\locb$ and of $\bullet$ on non-loop edges must be even, as the total number of $\locb$ and of $\circ$ on non-loop edges must be even too;

    \item  there cannot be a $\circ$-decoration to the left of a $\bullet$-decoration and vice versa. 
\end{enumerate}
We divide the admissible diagrams in five disjoint families, based on the displacement of the decorations on the edges (ALT, LP, RP, LRP and PZZ-diagrams). As one can guess, there will be a correspondence between this diagram classification and the heaps classification in Section \ref{fc}. Moreover we give a definition of the \textit{length of a diagram} $\ell(D)$ that depends on the shape of the edges and on the number of decorations on them. We use the length of a diagram for the inductive argument of our main result. The following picture summarizes the several structures introduced above: the main result is described by the last equality.

\begin{center}
\begin{tikzcd}
    T_{n+2}^{LR}(\Omega)
    \arrow[r, hook] &  \mathcal{P}_{n+2}^{LR}(\Omega) \arrow [r, twoheadrightarrow] & \PLRn \supset \DD = \dbD.
\end{tikzcd}
\end{center}

\section{Cut and paste operation}
In this section we present the major combinatorial technique we used to prove the faithfulness of the diagrammatic representation we just defined.

	An edge on the north face $e$ of $D\in \dbD$ joining two consecutive edges is called a \textit{simple edge} if either (a) $e$ is undecorated, or (b) $e=\{1,2\}$ and it is decorated by a single $\bullet$, or (c) $e=\{n+1, n+2\}$ and it is decorated by a single $\circ$.
	
\begin{figure}[h]
    \centering
    \includegraphics[scale=0.6]{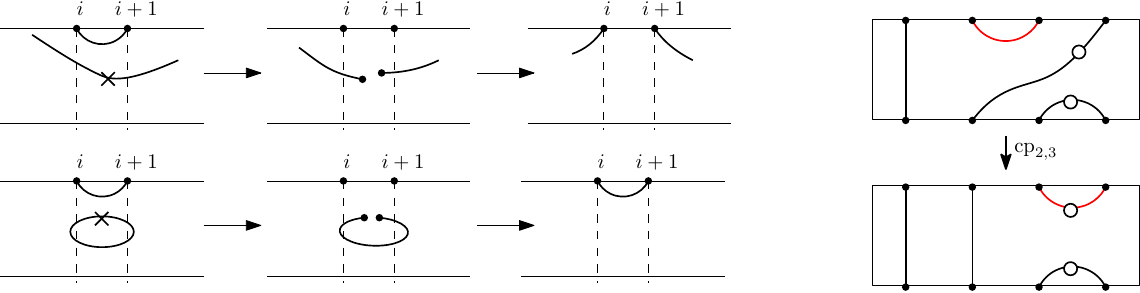}
    \caption{Cut and paste operation.}
    \label{fig:cutpaste}
\end{figure}

We consider a subset of simple edges, called \textit{suitable} edges, and to each of those, we assign a \textit{neighbor} edge which often is a $\locb$ or leaves one of the nodes adjacent to the suitable edge. Now we can introduce the following procedure, called \textit{cut and paste operation}, an example is in Figure \ref{fig:cutpaste}, right.

\begin{Definition}[Cut and paste operation]\label{definition:cutandpaste}
Let $D$ be an admissible diagram with a suitable edge $e=\{i, i+1\}$.
\begin{enumerate}
\item[(cp1)] Delete the simple edge $e$.
\item[(cp2)] Cut the neighbor of $e$ and join the two free endpoints of the cut edge to the nodes $i$ and $i+1$ as to obtain two new non intersecting edges or the edge $\{i,i+1\}$, see Figure~\ref{fig:cutpaste}, left. 
\item[(cp3)] If the neighbor edge of $e$ was decorated, then distribute its decorations on the new edges in an admissible way.
\end{enumerate}
\end{Definition} 

The importance of this procedure relies in the possibility of finding a unique diagram $D'$ from another diagram $D$ that factorizes as $D=D_iD'$ with the property that $D'$ has length equal to $\ell(D)-1$. Thanks to this inductive argument, we can prove the results stated in the next section. Moreover, this procedure provides an algorithm to factorize an admissible diagram into a product of simple diagrams.

\section{Main results}
Let $s_{i_1}s_{i_2}\cdots s_{i_k}$ be a reduced expression of $w\in \FC(\Dn)$ and define $D_w:=D_{i_1}D_{i_2}\cdots D_{i_k}$. Note that $D_w$ does not depend on the chosen reduced expression of $w$ since $w\in \FC(\Dn)$. Define the 
$\Zd$-algebra homomorphism $$\tilde{\theta}_D: \TLD \rightarrow \DD \qquad \mbox{ such that }\qquad \tilde{\theta}_D(b_i)=D_i$$ for all $i=0,\ldots,n+2$. 
Clearly $\tilde{\theta}_D$ is surjective and maps the monomial basis element $b_w:=b_{i_1}b_{i_2}\cdots b_{i_k}$ into the diagram $D_w$. Our goal is to show that this map is actually an algebra isomorphism. The proof is based on five steps that are summarized in the next theorem. Point 4) is by induction on the length of an admissible diagram and it uses the cut and paste algorithm. 

\begin{Theorem}\label{factorization}\

    \begin{enumerate}
        \item Let $w$ be a FC element of a certain type, then the image of       the basis element $b_w$ is a diagram of the analogous type and vice versa (for instance, $w$ has a heap of type (ALT) if and only if $\tilde{\theta}_D(b_w)=D_w$ is an ALT-diagram, and so on). 
\begin{figure}[h]
    \centering
    \includegraphics[scale = 0.7]{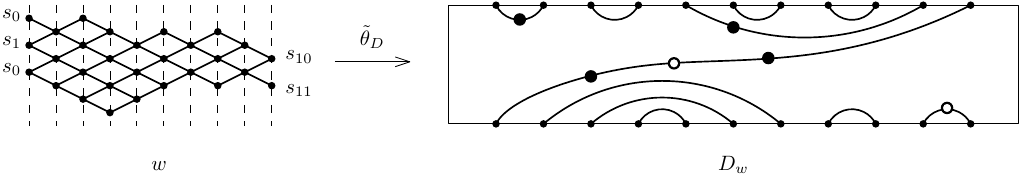}
    \label{fig:enter-label}
\end{figure}

        \item The lengths of $w$ and $D_w$ are equal.
        \item Every admissible diagram $D$ is of the form $D_w$ where $w\in \FC(\Dn)$.
        \item The admissible diagrams form a $\Zd$-algebra that coincides with $\DD$. Moreover, the set of admissible diagrams is a basis for $\DD$.
        \item The map $\tilde{\theta}_D$ is an algebra isomorphism.
    \end{enumerate}
	\end{Theorem}

\nocite{*}
\bibliographystyle{eptcs}
\bibliography{generic}

\begin{thebibliography}{10}
\providecommand{\bibitemdeclare}[2]{}
\providecommand{\surnamestart}{}
\providecommand{\surnameend}{}
\providecommand{\urlprefix}{Available at }
\providecommand{\url}[1]{\texttt{#1}}
\providecommand{\href}[2]{\texttt{#2}}
\providecommand{\urlalt}[2]{\href{#1}{#2}}
\providecommand{\doi}[1]{doi:\urlalt{https://doi.org/#1}{#1}}
\providecommand{\eprint}[1]{arXiv:\urlalt{https://arxiv.org/abs/#1}{#1}}
\providecommand{\bibinfo}[2]{#2}

\bibitemdeclare{article}{BJNFC}
\bibitem{BJNFC}
\bibinfo{author}{R.~\surnamestart Biagioli\surnameend},
  \bibinfo{author}{F.~\surnamestart Jouhet\surnameend} \&
  \bibinfo{author}{P.~\surnamestart Nadeau\surnameend} (\bibinfo{year}{2015}):
  \emph{\bibinfo{title}{Fully commutative elements in finite and affine
  {C}oxeter groups}}.
\newblock {\slshape \bibinfo{journal}{Monatsh. Math.}}
  \bibinfo{volume}{178}(\bibinfo{number}{1}), pp. \bibinfo{pages}{1--37},
  \doi{10.1007/s00605-014-0674-7}.
\newblock \eprint{1402.2166}.

\bibitemdeclare{article}{ErnstDiagramI}
\bibitem{ErnstDiagramI}
\bibinfo{author}{D.C. \surnamestart Ernst\surnameend} (\bibinfo{year}{2012}):
  \emph{\bibinfo{title}{Diagram calculus for a type affine {$C$}
  {T}emperley--{L}ieb algebra, {I}}}.
\newblock {\slshape \bibinfo{journal}{J. Pure Appl. Algebra}}
  \bibinfo{volume}{216}(\bibinfo{number}{11}), pp. \bibinfo{pages}{2467--2488},
  \doi{10.1016/j.jpaa.2012.03.013}.
\newblock \eprint{0910.0925}.

\bibitemdeclare{article}{ErnstDiagramII}
\bibitem{ErnstDiagramII}
\bibinfo{author}{D.C. \surnamestart Ernst\surnameend} (\bibinfo{year}{2018}):
  \emph{\bibinfo{title}{Diagram calculus for a type affine {$C$}
  {T}emperley--{L}ieb algebra, {II}}}.
\newblock {\slshape \bibinfo{journal}{J. Pure Appl. Algebra}}
  \bibinfo{volume}{222}(\bibinfo{number}{12}), pp. \bibinfo{pages}{3795--3830},
  \doi{10.1016/j.jpaa.2018.02.008}.
\newblock \eprint{1101.4215}.

\bibitemdeclare{article}{FanGreenAffine}
\bibitem{FanGreenAffine}
\bibinfo{author}{C.K. \surnamestart Fan\surnameend} \& \bibinfo{author}{R.M.
  \surnamestart Green\surnameend} (\bibinfo{year}{1999}):
  \emph{\bibinfo{title}{On the affine {T}emperley--{L}ieb algebras}}.
\newblock {\slshape \bibinfo{journal}{J. London Math. Soc.}}
  \bibinfo{volume}{60}(\bibinfo{number}{2}), pp. \bibinfo{pages}{366--380},
  \doi{10.1112/S0024610799007796}.
\newblock \eprint{q-alg/9706003}.

\bibitemdeclare{phdthesis}{Graham}
\bibitem{Graham}
\bibinfo{author}{J.J. \surnamestart Graham\surnameend} (\bibinfo{year}{1995}):
  \emph{\bibinfo{title}{Modular Representations of Hecke Algebras and Related
  Algebras}}.
\newblock Ph.D. thesis, \bibinfo{school}{University of Sydney}.

\bibitemdeclare{article}{Greencellular}
\bibitem{Greencellular}
\bibinfo{author}{R.M. \surnamestart Green\surnameend} (\bibinfo{year}{1998}):
  \emph{\bibinfo{title}{Cellular algebras arising from {H}ecke algebras of type
  {$H_n$}}}.
\newblock {\slshape \bibinfo{journal}{Math. Z.}}
  \bibinfo{volume}{229}(\bibinfo{number}{2}), pp. \bibinfo{pages}{365--383},
  \doi{10.1007/PL00004661}.
\newblock \eprint{q-alg/9712019}.

\bibitemdeclare{article}{Greengeneral}
\bibitem{Greengeneral}
\bibinfo{author}{R.M. \surnamestart Green\surnameend} (\bibinfo{year}{1998}):
  \emph{\bibinfo{title}{Generalized {T}emperley--{L}ieb algebras and decorated
  tangles}}.
\newblock {\slshape \bibinfo{journal}{J. Knot Theory Ramifications}}
  \bibinfo{volume}{7}(\bibinfo{number}{2}), pp. \bibinfo{pages}{155--171},
  \doi{10.1142/S0218216598000103}.
\newblock \eprint{q-alg/9712018}.

\bibitemdeclare{article}{GreenStar}
\bibitem{GreenStar}
\bibinfo{author}{R.M. \surnamestart Green\surnameend} (\bibinfo{year}{2006}):
  \emph{\bibinfo{title}{Star reducible {C}oxeter groups}}.
\newblock {\slshape \bibinfo{journal}{Glasg. Math. J.}}
  \bibinfo{volume}{48}(\bibinfo{number}{3}), pp. \bibinfo{pages}{583--609},
  \doi{10.1017/S0017089506003211}.
\newblock \eprint{math/0509363}.

\bibitemdeclare{article}{GreenTLE}
\bibitem{GreenTLE}
\bibinfo{author}{R.M. \surnamestart Green\surnameend} (\bibinfo{year}{2009}):
  \emph{\bibinfo{title}{On the {M}arkov trace for {T}emperley--{L}ieb algebras
  of type {$E_n$}}}.
\newblock {\slshape \bibinfo{journal}{J. Knot Theory Ramifications}}
  \bibinfo{volume}{18}(\bibinfo{number}{2}), pp. \bibinfo{pages}{237--264},
  \doi{10.1142/S0218216509006872}.
\newblock \eprint{0704.0283}.

\bibitemdeclare{article}{Joneshecke}
\bibitem{Joneshecke}
\bibinfo{author}{V.F.R. \surnamestart Jones\surnameend} (\bibinfo{year}{1987}):
  \emph{\bibinfo{title}{Hecke algebra representations of braid groups and link
  polynomials}}.
\newblock {\slshape \bibinfo{journal}{Ann. of Math.}}
  \bibinfo{volume}{126}(\bibinfo{number}{2}), pp. \bibinfo{pages}{335--388},
  \doi{10.2307/1971403}.

\bibitemdeclare{article}{Kauffman}
\bibitem{Kauffman}
\bibinfo{author}{L.H. \surnamestart Kauffman\surnameend}
  (\bibinfo{year}{1987}): \emph{\bibinfo{title}{State models and the {J}ones
  polynomial}}.
\newblock {\slshape \bibinfo{journal}{Topology}}
  \bibinfo{volume}{26}(\bibinfo{number}{3}), pp. \bibinfo{pages}{395--407},
  \doi{10.1016/0040-9383(87)90009-7}.

\bibitemdeclare{article}{Penrose}
\bibitem{Penrose}
\bibinfo{author}{R.~\surnamestart Penrose\surnameend} (\bibinfo{year}{1971}):
  \emph{\bibinfo{title}{Angular momentum: an approach to combinatorial
  spacetime}}.
\newblock {\slshape \bibinfo{journal}{Quantum theory and beyond}}
  \bibinfo{volume}{151}, pp. \bibinfo{pages}{395--407}.

\bibitemdeclare{article}{TemperleyLieb}
\bibitem{TemperleyLieb}
\bibinfo{author}{H.N.V. \surnamestart Temperley\surnameend} \&
  \bibinfo{author}{E.H. \surnamestart Lieb\surnameend} (\bibinfo{year}{1971}):
  \emph{\bibinfo{title}{Relations between the ``percolation'' and ``colouring''
  problem and other graph-theoretical problems associated with regular planar
  lattices: some exact results for the ``percolation'' problem}}.
\newblock {\slshape \bibinfo{journal}{Proc. Roy. Soc. London Ser. A}}
  \bibinfo{volume}{322}(\bibinfo{number}{1549}), pp. \bibinfo{pages}{251--280},
  \doi{10.1098/rspa.1971.0067}.

\end{thebibliography}
\end{document}